\documentclass[12pt]{article}
 \usepackage{a4, latexsym, amsfonts}
 \newtheorem{theorem}{Theorem}[section]

 \newtheorem{conjecture}{Conjecture}[section]

\newcommand{\supp}{{\rm supp}}

 \newenvironment{proof}{\trivlist
      \item[\hskip\labelsep
      {\itshape Proof.}]\normalfont}
      {\hspace*{\fill}$\Box$\endtrivlist}

\begin{document}

\title{Davenport constant with weights}
\author{  Pingzhi Yuan \\
{\small School of Mathematics, South China Normal University
  , Guangzhou 510631, P.R.CHINA}\\
 {\small e-mail  mcsypz@mail.sysu.edu.cn}\\ Xiangneng  Zeng\\
\small{ Department of Mathematics, Sun Yat-Sen University, Guangzhou
510275, P.R.CHINA}}

\date{}
\maketitle
 \edef \tmp {\the \catcode`@}
   \catcode`@=11
   \def \@thefnmark {}

    \@footnotetext {  Supported by the Guangdong Provincial Natural Science Foundation (No. 8151027501000114) and  NSF of China (No. 10571180).}

\begin{abstract}
For the cyclic group $G=\mathbb{Z}/n\mathbb{Z}$ and any non-empty
$A\in\mathbb{Z}$.    We define the Davenport constant of $G$ with
weight $A$, denoted by $D_A(n)$, to be the least natural number $k$
such that for any sequence $(x_1, \cdots, x_k)$ with $x_i\in G$,
there exists a non-empty subsequence $(x_{j_1},\, \cdots,\, x_{j_l}
)$ and $a_1, \, \cdots,\,  a_l\in A$ such that $\sum_{i=1}^l
a_ix_{j_i} = 0$. Similarly, we define the constant $E_A(n)$ to be
the least $t\in\mathbb{N}$ such that for all sequences $(x_1, \,
\cdots,\, x_t )$ with $x_i \in G$, there exist indices $j_1,\,
\cdots, \,j_n\in\mathbb{N}, 1\leq j_1 <\cdots < j_n\leq t$, and
$\vartheta_1,\, \cdots,\,\vartheta_n\in A$ with $\sum^{n}_{ i=1}
\vartheta_ix_{j_i} = 0$. In the present paper, we show that
$E_A(n)=D_A(n)+n-1$. This solve the problem raised by Adhikari and
Rath \cite{ar06}, Adhikari and Chen \cite{ac08}, Thangadurai
\cite{th07} and  Griffiths \cite{gr08}.


\end{abstract}

\section{Introduction}

 For
an abelian group $G$, the Davenport constant $D(G)$ is defined to be
the smallest natural number $k$ such that any sequence of $k$
elements in $G$ has a non-empty subsequence whose sum is zero (the
identity element). Another interesting constant $E(G)$ is defined to
be the smallest natural number $k$ such that any sequence of $k$
elements in $G$ has a subsequence of length $|G|$ whose sum is zero.

The following result due to Gao \cite{ga96} connects these two
invariants.

\begin{theorem} If $G$ is a finite abelian group of order $n$, then $E(G)
= D(G) + n -1$.\end{theorem}

For a finite abelian group $G$ and any non-empty $A\in\mathbb{Z}$,
Adhikari and Chen \cite{ac08}  defined the Davenport constant of $G$
with weight $A$, denoted by $D_A(G)$, to be the least natural number
$k$ such that for any sequence $(x_1, \cdots, x_k)$ with $x_i\in G$,
there exists a non-empty subsequence $(x_{j_1},\, \cdots,\, x_{j_l}
)$ and $a_1, \, \cdots,\,  a_l\in A$ such that $\sum_{i=1}^l
a_ix_{j_i} = 0$. Clearly, if $G$ is of order $n$, it is equivalent
to consider $A$ to be a non-empty subset of $\{0, 1, \cdots , n
-1\}$ and cases with $0\in A$ are trivial.

Similarly, for any such set $A$, for a finite abelian group $G$ of
order $n$, the constant $E_A(G)$ is defined to be the least
$t\in\mathbb{N}$ such that for all sequences $(x_1, \, \cdots,\, x_t
)$ with $x_i \in G$, there exist indices $j_1,\, \cdots,
\,j_n\in\mathbb{N}, 1\leq j_1 <\cdots < j_n\leq t$, and
$\vartheta_1,\, \cdots,\,\vartheta_n\in A$ with $\sum^{n}_{ i=1}
\vartheta_ix_{j_i} = 0$.

For the group $G = \mathbb{Z}/n\mathbb{Z}$, we write $E_A(n)$ and
$D_A(n)$ respectively for $E_A(G)$ and $D_A(G)$. In the cases
$A=\{1\}, \, \{-1, \, 1\}, \mathbb{Z}_n^\star$  or $A=(a_1, \cdots,
a_r)$ with $\gcd(a_2-a_1, \cdots, a_r-a_1, n)=1$ or $n=p$ is a
prime, it is proved that $E_A(n)= D_A(n)+n-1$. The following
conjecture has been raised by Adhikari and Rath \cite{ar06},
Adhikari and Chen \cite{ac08}, Thangadurai \cite{th07} and Griffiths
\cite{gr08}, they seems to believe that Conjecture 1.1 is true and
have proved it in some special cases.

\begin{conjecture}
For any non-empty set $A\in\mathbb{Z}$, $E_A(n)=D_A(n)+n-1$.
\end{conjecture}

The main purpose of the present paper is to prove Conjecture 1.1. By
using the main theorem of Devos, Goddyn and Mohar \cite{dgm08} and a
recently proved theorem of the authors \cite{yz08}, we shall prove
the following theorem
\begin{theorem}

For any non-empty set $A\in\mathbb{Z}$, $E_A(n)=D_A(n)+n-1$.
\end{theorem}

Throughout this paper, let $G$ be an additive finite abelian group.
$\mathcal{F}(G)$ denotes the free abelian monoid with basis $G$, the
elements of which are called $sequences$ (in $G$). A sequence of not
necessarily distinct elements from $G$ will be written in the form
$S=g_1\, \cdots \, g_k=\prod_{i=1}^k g_i=\prod_{g\in G}g^{\mathsf
v_g(S)}\in\mathcal{F}(G)$, where $\mathsf v_g(S) \geq 0$ is called
the $multiplicity$ of $g$ in $S$. We call $|S|= k$ the $length$ of
$S$, $\mathsf h(S)=\max\{\mathsf v_g(S)|g\in G\}\in[0, |S|]$ the
maximum of the multiplicities of $S$,   $\supp(S) = \{g\in G: \,
\mathsf v_g(S)>0\}$  the $support$ of $S$. For every $g\in G$ we set
$g+S=(g+g_1) \cdots(g+g_k)$.

We say that $S$ contains some $g\in G$ if $\mathsf v_g(S)\geq1$ and
a sequence $T\in\mathcal{F}(G)$ is a $subsequence$ of $S$ if
$\mathsf v_g(T) \leq \mathsf v_g(S)$ for every $g\in  G$, denoted by
$T|S$.   Furthermore, by $\sigma(S)$ we denote the sum of $S$, (i.e.
$\sigma(S) = \sum_{ i=1}^k g_i = \sum_{g\in G}\mathsf v_g(S)g\in
G$). For every $k\in\{1, 2, \cdots, \, |S|\}$, let
$\sum_k(S)=\{g_{i_1}+\cdots+g_{i_k}|1\leq i_1<\cdots<i_k\leq|S|\},
\, \sum_{\leq k}(S)=\cup_{i=1}^k\sum_{i}(S)$, and let
$\sum(S)=\sum_{\leq |S|}(S)$.

Let $S$ be a sequence in $G$. We call $S$ a $zero-sum$ $sequence$ if
$\sigma(S) =
 0$.

Also, we follow the same terminologies and notations as in the
survey article \cite{gg06} or in the book \cite{gh06}.

\section{Lemmas}
First, we need  a result on the sum of $l$ finite subsets of $G$. If
${\bf A}=(A_1, A_2, \cdots, A_m)$ is a sequence of finite subsets of
$G$, and $l\leq m$, we define
$$\sum_l({\bf A})=\{a_{i_1}+\cdots+a_{i_l}: 1\leq i_1<\cdots<i_l\leq
m \, \mbox{ and}\, a_{i_j}\in A_{i_j} \, \mbox{ for every}\, 1\leq
j\leq l\}.$$ So $\sum_l({\bf A})$ is the set of all elements which
can be represented as a sum of $l$ terms from distinct members of
${\bf A}$. The following is the main result of Devos, Goddyn and
Mohar \cite{dgm08}.

{\bf Theorem DGM} {\it Let ${\bf A}=(A_1, A_2, \cdots, A_m)$ be a
sequence of finite subsets of $G$, let $l\leq m$, and let
$H=stab(\sum_l({\bf A}))$.  If $\sum_l({\bf A})$ is nonempty, then
$$|\sum_l({\bf A})|\geq |H|(1-l+\sum_{Q\in G/H}\min\{l, |\{i\in\{1,
\cdots, m\}:A_i\cap Q\neq\emptyset\}|\}).$$}

We still need the following new result on Davenport's constant
\cite{yz08}.

{\bf Theorem YZ} {\it Let $G$ be a finite abelian group of order $n$
and Davenport constant $D(G)$. Let $S=0^{\mathsf h(S)}\prod_{g\in
G}g^{\mathsf v_g(S)}\in \mathcal{F}(G)$ be  a sequence with a
maximal multiplicity $\mathsf h(S)$ attained by $0$ and  $|S|=t\geq
n+D(G)-1$. Then there exists a subsequence $S_1$ of $S$ with length
$|S_1|\geq t+1-D(G)$ and $0\in \sum_k(S_1)$ for every $1\leq k\leq
|S_1|$. In particular, for every sequence $S$ in $G$ with length
$|S|\geq n+D(G)-1$, we have
$$0\in \sum_{km}(S), \,\mbox{ for every } \quad 1\leq k\leq (|S|+1-D(G))/m,$$
where $m$ is the exponent of $G$.}
\section{Proof of Theorem 1.2}

\begin{proof} The proof of $E_A(n)\geq D_A(n)+n-1$ is easy, so it is
sufficient to prove the reverse inequality.

For any non-empty set $A=\{a_1, \cdots a_r\}\subset\mathbb{Z}$ and a
cyclic group $G=\mathbb{Z}/n\mathbb{Z}$, let $t=D_A(n)+n-1$ and
$S=x_1\cdots x_t$ is any sequence in $G$ with length
$|S|=t=D_A(n)+n-1$. Put
$$A_i=Ax_i=\{a_1x_i, \cdots, a_rx_i\}  \mbox{ for}\, i=1, \cdots,
t$$ and ${\bf A}=(A_1, \, \cdots, \, A_t)$. It suffices to prove
that $0\in\sum_n({\bf A})$.

 We shall assume (for a contradiction) that
the theorem is false and choose a counterexample $(A, G, \, S)$ so
that $n=|G|$ is minimum, where $G$ is a cyclic group of order $n$,
$A$ is a finite subset of $\mathbb{Z}$ and  $S=x_1\cdots x_t$ is a
sequence in $G$ such that
$$0\not\in\sum_n({\bf A}).$$

Next we will show that our assumptions imply $H=stab(\sum_n({\bf
A}))=\{0\}$. Suppose (for a contradiction) that $H=stab(\sum_n({\bf
A}))\neq\{0\}$ and let $\varphi: \, G \longrightarrow G/H$ denote
the canonical homomorphism and $\varphi(x_i)$ the image of $x_i$ for
$1\leq i \leq t$. Let ${\bf A_\varphi}= (\varphi(A_1), \cdots,
\varphi(A_t))$. By our assumption for the minimal of $|G|$, the
theorem holds for $(A, \varphi(G), \, \varphi(S))$. Since
$n>|\varphi(G)|+D_A(\varphi(G))-1$, $|\varphi(G)||n$ and $D_A(G)\geq
D_A(\varphi(G))$, repeated applying the theorem  to the sequence
$\varphi(S)=\varphi(x_1), \cdots, \varphi(x_t)$ we have
$$\varphi(0)=\varphi(H)\in\sum_n({\bf A_\varphi}),$$
thus $0\in H\subset\sum_n({\bf A})$. This  contradiction implies
that $H=stab(\sum_n({\bf A}))=\{0\}$.

If there is an element $a\in G$ such that $|\{j\in\{1, \cdots, t\}:
a\in A_j\}|\geq n$, then $0\in\sum_n({\bf A})$, a contradiction.
Therefore we may assume that for every $a\in G$, $|\{j\in\{1,
\cdots, t\}: a\in A_j\}|\leq n$. Let  $r$ be the number of
$i\in\{1,\cdots, t\}$ with $|A_i|=1$, by Theorem DGM and the
assumptions, we have

$$n-1\geq \sum_n({\bf A})\geq 1-n +\sum_{a\in G}\min\{n,|\{j\in\{1, \cdots, t\}: a\in A_j\}|\}$$
$$=1-n+\sum_{i=1}^t|A_i|\geq 1-n+2(n+D_A(G)-1-r)+r.$$
 It
follows that \begin{equation}r\geq 2D_A(G).\end{equation}

Without loss of generality, we may assume that $x_1, \cdots, x_r$
are all the elements in $\{x_1, \, \cdots, \, x_t\}$ such that
$|A_i|=1$, and $x_1$ is the element in $\{x_1, \, \cdots, \, x_r\}$
such that $a_1x_1$ attains the maximal multiplicity in the sequence
$S_1=(a_1x_1)\cdots(a_1x_r)$. Observe that $\sum_n({\bf A})=
\sum_n(A(x_1-x_u), \, \cdots, \, A(x_t-x_u))$ for every $1\leq u\leq
r$. Therefore without loss of generality we may assume that
$a_1x_1=0$ and $\mathsf v_0(S_1)=\mathsf h(S_1)$ for the sequence
\begin{equation}
S_1=(a_1x_1)\cdots(a_1x_r)=0^{\mathsf h(S_1)}(a_1x_{\mathsf
h(S_1)+1})\cdots(a_1x_r). \end{equation} Let $H_1=<x_1, \cdots,
x_r>$ be the group generated by  $x_1, \cdots, x_r$, $H=a_1H_1$. We
have the following claim.

{\bf Claim:} $D_A(G)\geq D_A(H_1)\geq D(H)=|H|$.

The last equality of the Claim follows from the fact that $H$ is a
subgroup of the cyclic group $G$. The first inequality  in the Claim
is obvious, so we only need to prove that $D_{A}(H_1)\geq D(H)$.
Suppose that  $W=y_1\cdots y_{D(H)-1}$ is a zero-sum free sequence
in $H$. Since $H=a_1H_1$, we have $ y_i=a_1w_i, w_i\in H_1, i=1,
\cdots r$. Further, it is easy to see that $ Aw_i=a_1w_i, \, i=1,
\cdots, r$ by the definition of $H_1$, so $w_1\cdots w_{D(H)-1}$ is
a zero-sum free sequence in $H_1$ with respect to the weight $A$,
thus $D_A(H_1)\geq D(H)$ and the Claim follows.

By the Claim, (1), (2) and Theorem YZ, $S_1$ has a subsequence $S_2$
of length $|S_2|=s\geq r+1-|H|$ such that $0\in\sum_l(S_2)$ for
every $1\leq l\leq s$. Without loss of generality, we may assume
that $S_2=(a_1x_1)\cdots(a_1x_s)$.

If $s\geq n$ then $0\in\sum_n(S_2)\subset\sum_n({\bf A})$, we are
done.

If $s<n$, then $|x_{s+1}\cdots x_t|=t-s=n-1+D_A(G)-s\geq D_A(G)$.
Repeated using  the definition of $D_A(G)$, there exists an integer
$v$ such that $v\leq n, \, t-s-v\leq D_A(G)-1$ and
$$0\in\sum_v((A_{s+1}, \cdots A_t)).$$
Since $\sum_l(S_2)=\sum_l((A_1, \cdots, A_s))$ and $0\in\sum_l((A_1,
\cdots, A_s))$ for every $1\leq l\leq s$, we have
$$0\in\sum_{v+k}({\bf A}) \quad \mbox{ for every } 0\leq k\leq s.$$
Therefore $0\in\sum_n({\bf A})$ since $v+s\geq t+1-G_A(G)\geq n$.
This completes the proof of the theorem.

\end{proof}

\end{document}